\documentclass[12pt]{amsart}
\usepackage{amsmath,amssymb,amsfonts,amsthm,amstext,verbatim}
\usepackage{url,enumerate,color,graphicx,stmaryrd,psfrag,ulem} 
\newcommand\seprule{\noalign{\vskip0.5ex}\hline\noalign{\vskip0.5ex}}

\newcounter{mycount}

\begin{document}
\title[Power-Weighted Variant of the Cambridge Compromise]{A Power-Weighted Variant of the\\EU27 Cambridge Compromise}

\author[Grimmett]{G.\ R.\ Grimmett}
\address{Statistical Laboratory, Centre for
Mathematical Sciences, Cambridge University, Wilberforce Road,
Cambridge CB3 0WB, UK} \email{g.r.grimmett@statslab.cam.ac.uk}
\urladdr{http://www.statslab.cam.ac.uk/$\sim$grg/}

\author[Oelbermann]{K.-F.\ Oelbermann}
\author[Pukelsheim]{F.\ Pukelsheim}
\address{Institut f\"ur Mathematik, Universit\"at Augsburg, 
D-86135 Augsburg, Germany}
\email{\{Oelbermann, Pukelsheim\}@math.uni-augsburg.de}
\urladdr{http://www.math.uni-augsburg.de/stochastik/pukelsheim/}
\urladdr{http://www.math.uni-augsburg.de/stochastik/oelbermann/}

\begin{abstract}
The Cambridge Compromise composition of the European Parliament
allocates five base seats to each Member State's citizenry, and apportions
the remaining seats proportionately to population figures using the divisor
method with rounding upwards and observing a 96 seat capping.  
The power-weighted variant avoids the capping step,  proceeding
instead  by a progressive non-linear downweighting of
the population figures until the largest State is allocated exactly 96 seats.
The pertinent calculations of the variant are described, and its
relative constitutional merits are discussed.
\end{abstract}

\date{4 August 2011}  %\today}

\subjclass[2010]{}

\maketitle

\section{Introduction}\label{sec:1}

The Cambridge Compromise ``5+Upw'' is a mathematical formula to
apportion the 751 seats of the European Parliament (EP) among the 27
Member States, \cite{grim}.  It operates in two steps.  In the first
step, five {\it base seats} are allocated to the citizenry of each
Member State.  In the second step, the $616\ (= 751 - 27 \cdot 5)$
{\it remaining seats} are apportioned in proportion to the States'
population figures\footnote{We find it useful to sort populations from
  largest to smallest.  Thus $p_1$ refers to the largest Member State
  (at present, Germany), and $p_{27}$ to the smallest (at present,
  Malta).}  $p_1, \ldots, p_{27}$.  The divisor method with rounding
upwards secures the minimum threshold of six seats per Member
State. Seat contingents are capped at 96 if need be, staying in line
with Article 14 II of the Treaty of Lisbon that decrees that no Member
State shall be allocated less than 6 or more than 96 seats,
\cite{teu}.

This note develops a {\it power-weighted variant} ``5+Pwr+Upw'' of the
Cambridge Compromise, in the spirit of \cite{arndt,theil}.  The
variant refers the apportionment of the 616 remaining seats to
power-weighted population indices $p_i^E$, rather than to the original
population figures $p_i$.  When necessary, an exponent $E<1$ is
calculated in order that the contingent of the largest Member State be
reduced to exactly 96 seats.

It is shown in Section~\ref{sec:2} that, for a given apportionment
vector \linebreak $x = (x_1, \ldots, x_{27})$, the exponent~$E$ ranges
across an interval $[E_{\min}(x), E_{\max}(x)]$ without alteration to
$x$.  We describe in Section~\ref{sec:3} how to determine an exponent
$E$ that secures 96 seats for the largest State.  Section~\ref{sec:4}
contains a discussion of the merits of the variant from the viewpoint
of general electoral principles, and of the requirements of the
Union's primary law.

The power-weighted variant stays closer to the {\it status quo}
composition of the EP than do (at least) many of the other methods reviewed in
\cite{grim}. It affords a possible transition from the present
composition to the eventual Cambridge Compromise ``5+Upw'' across an
enlarged Union, and it satisfies the amended definition of degressive
proportionality of \cite{grim}.  On the other hand, the use of
power-weighted population indices entails a lack of transparency, and
constitutes a breach of the principle of equal voting power.

\section{Exponent-ranges}\label{sec:2}

Suppose the {\it seat vector} $x$ originates from population indices
that are weighted with some given exponent $E > 0$.  Since the
exponent is a continuous variable while the seat contingents are
integer-valued, there exists a range of values of~$E$ across which the
seat vector $x$ is constant.  In this section, we investigate the
induced {\it exponent-range} $[E_{\min}(x), E_{\max}(x)]$.  The allocation of base
seats plays no role in this investigation, and will therefore be
neglected in this section.  Similarly, we shall take no account of the
96-seat cap in this section.

The seats to be apportioned are the 616 remaining seats.  Let $y =
y(E)$ be a seat apportionment vector generated with a given exponent
$E > 0$.  The components $y_i$ are obtained from the population
indices $p_i^E$ via scaling with a common divisor $D > 0$ and rounding
the resulting quotients upwards, where the divisor is chosen so that
the sum of the components exhausts the seats available:
   $$
   y_i = \left \lceil { p_i^E \over D} \right \rceil \quad \hbox{ for
   }i = 1, \ldots, 27, \qquad \sum_{i = 1}^{27} y_i = 616.
   $$
For a given $y$, the exponent-range $[E_{\min}(y), E_{\max}(y)]$ is
calculated as follows. 

We need consider only seat vectors $y$ with $y_i \ge 1$ for all $i$.
Consider two States indexed $i$ and $j$.  There exists a unique {\it
  critical exponent} $E(i,j)$, depending on $y$, such that, when
$E=E(i,j)$, the allocations $(y_i, y_j)$ and $(y_i+1, y_j-1)$ are
equally justified. Only situations with $y_j\geq 2$ are relevant in
order to make sure that $y_j-1\geq 1$ plus five base seats achieve the
six seat minimum required by the Treaty of Lisbon.

The tie entails the identities
   $$ 
   { p_i^{E(i,j)} \over D} = y_i \quad \hbox{ and } \quad
   { p_j^{E(i,j)} \over D} = y_j-1.
   $$
In view of the rule of rounding upwards, this indicates that State
$i$ may be allocated $y_i$ or $y_i+1$ seats, and that State $j$ is
eligible for $y_j-1$ or $y_j$ seats.  The identities allow us to
eliminate the divisor and to solve thus for the critical exponent,
   $$
   E(i,j) = { \displaystyle\log{(y_i/( y_j-1))} 
        \over \displaystyle\log{(p_i / p_j)} }.
   $$ 
As $E$ increases, the seat vector $y$ is unchanged until $E$ attains the
value
   $$
   E_{\max}(y) = \min \left\{  E(i,j) : p_i > p_j, \, y_j\geq 2\right\}. 
   $$
By a similar argument as $E$ decreases, the lower boundary point of
the exponent-range is
   $$
   E_{\min}(y) = \max \left\{  E(i,j) : p_i < p_j, \, y_j\geq 2 \right\}.
   $$
The exponent-range for $y$ is found to be $[E_{\min}(y),
E_{\max}(y)]$.  This range necessarily contains the given exponent $E$ used
to generate the seat vector $y$. 

For example, consider the case of unweighted population figures,
$E=1$.  The seat numbers $x_i(1) = 5 + y_i(1)$ are exhibited in the
penultimate column of Table~\ref{table1}.  It transpires that $y(E) =
y(1)$ is unchanged as $E$ ranges from $E_{\min}\bigl(y(1)\bigr) =
0.9956$ to $E_{\max}\bigl(y(1)\bigr) = 1.0010$.  The \lq\lq nicest"
value in this exponent-range is evidently $E=1$.

For a general compact interval $[E_{\min}, E_{\max}]$, we may pick a
``nice'' representative exponent~$E$ by rounding the midpoint of the
interval to as few significant digits as the interval permits.  The
smallest and greatest intervals are half-open, however.  In these two
intervals, we may choose a representative exponent near the closed
boundary, in order to avoid numerical difficulties.

The transfer argument proves also that, if $E < E'$, the seat vector
$y(E)$ is majorized\footnote{The property of majorization amounts to
  the fact that, for any $k$, the aggregate number of seats of the $k$
  largest States is nondecreasing.  As a consequence, the aggregate
  number of seats of the $27-k$ smallest States is nonincreasing.  See
  \cite{marshall} for more details.}  by the seat vector $y(E')$.  In
particular, the largest State $i = 1$ has seat numbers that are
nondecreasing in $E$, that is $E < E' \Rightarrow y_1(E) \le y_1(E')$.
We have that $y_1(0.01) = 28 - 5 = 23$ and $y_1(27.5) = 595 - 5 =
590$, see Table~\ref{table1}.  Hence, the equation $y_1(E) = 91$ is
solvable for $E$.  We consider the determination of a solution in the
next section.

\section{Choice of initial exponent}\label{sec:3}

We may determine an exponent $E$ with $y_1(E) = 91$ in the following algorithmic
manner.  As starting point, we may expect the ideal share of the
largest State to be close to its target contingent of 91
seats,
   $$
   { p_1^E \over \sum_{i =1}^{ 27} p_i^E} \cdot 616 = 91.
   $$ 
   This equation may be solved numerically, and yields (with the
   1.1.2011 Eurostat population figures) an initial exponent $E_{\rm
     init} = 0.8888$.  The largest State $i = 1$ misses its target by
   a single seat, $y_1(0.8888) = 90$.  The induced exponent-range
   turns out to be $[0.8884, 0.8977]$, with \lq\lq nice" exponent
   $0.893$.  The machinery of Section~\ref{sec:2} yields the next seat
   vector $y(0.899)$, with exponent-range $[0.8978, 8998]$, which
   still fails to allocate 91 seats to the largest State.  One further
   iteration yields an exponent $0.9 \in [0.8999, 0.9035]$ with which
   the discrepancy vanishes.  See Table~\ref{table1}.

   An improved initialization procedure is available.  The divisor
   method with rounding upwards is known to be biased, in that it has
   a tendency on average to favour smaller States at the expense of
   larger States.  The seat-bias of the largest State may be
   approximated by the formula
   $$
   - {1 \over 2} 
   \left\{ \left( \sum_{n =1}^{ 27} {1 \over n} \right) - 1 \right\}
   = - 1.4457,
   $$
   see \cite{heinrich}.  With this term included, the above
   initialization equation becomes
   $$
   { p_1^E \over \sum_{i =1}^{ 27} p_i^E} \cdot 616 = 91 + 1.4457,
   $$
   with numerical solution $E_{\rm init} = 0.9055$, and hence the
   \lq\lq nice" exponent $0.91$.  The ensuing seat vector $y(0.91)$
   achieves the target for the largest State, $y_1(0.91) = 91$, with
   exponent-range $[0.9036, 0.9109]$. Thus $y(0.91)$ differs from
   $y(0.9)$.

   We return now to the original setting in which each State receives
   in addition five base seats, and we write $x_i = 5 + y_i$.  If an
   exponent lies below 0.8999 or above 0.9109, the corresponding
   allocation fails to allocate to the largest State its target
   contingent of 96 seats.  The exponent-range $[0.8999, 0.9035]$
   gives rise to the feasible seat vector $x(0.9)$, and the subsequent
   exponent-range $[0.9036, 0.9109]$ yields the apportionment
   $x(0.91)$.  Thus we obtain two apportionment vectors, $x(0.9)$ and
   $x(0.91)$, each of which allocates 96 seats to the largest Member
   State.  See\footnote{Exponent-ranges for $E$ and applicable
     divisors $D$ are:\newline {\centering
\begin{tabular*}{\textwidth}{lllllllll}
  $E_{\min}$ &0&0.8884&0.8978&0.8999&0.9036&0.9110&0.9956&27.2202\\
  Exponent&0.01&0.893&0.899&0.9&0.91&0.912&1&28\\
  $E_{\max}$&0.0123&0.8977&0.8998&0.9035&0.9109&0.9125&1.0010&$\infty$\\ 
  Divisor&0.0526&121\,400&144\,400&146\,960&174\,600&180\,800&830\,000%
  &$6.12\cdot10^{218}$\\
\end{tabular*}}}
Table~\ref{table1}.

We close this section with some comments based on the current and
foreseeable population profile of the European Union.

1. There can exist several values of $E$ for which the seat vector $x(E)$ 
satisfies $x_1(E) = 96$.  Of these, the
smallest such exponent leads to the composition closest to the {\it status quo} composition and,
presumably, risks the greatest acclaim of the incumbent Parliament.
Similarly, the composition with the largest exponent comes
closest to the principled approach of the Cambridge Compromise.  See
Table~\ref{table2}.

2. A composition based on an exponent satisfying $E \le 1$
automatically satisfies the revised definition of degressive
proportionality of \cite{grim}.  Conversely, when $E>1$, the
representation might possibly be called ``progressive".

3. The power-weighted variant is a smoother allocation in situations
where the Cambridge Compromise ``5+Upw'' hits the upper cap.  Its
implementation is thus in two steps. First, calculate the CamCom
apportionment.  If this caps the largest State at 96, then calculate
the power-weighted variant.

4. As new States accede, the exponents $E$ that achieve $x_1(E) = 96$ 
approach unity.  When unity is reached, $E = 1$, the downweighting of 
population figures becomes neutral, and the ensuing apportionment
is that of the Cambridge Compromise.

5. The identity 5+Upw = 5.5+Std = 6+Dwn (\cite{zachariasen}) is
invalid for {\it  non-linear} variants of the Cambridge Compromise including
that considered here.

6. Serious issues of transparency and interpretation arise in the use
of power-weighted population indices.  For example, in the composition
with exponent $E = 0.9$, the number of (non-base) seats allocated to
Italy is in proportion to a population index of $10\,058\,816.8$
power-weighted ``apportionment units'', rather than to the Italian
population $60\,340\,328$.  How should this be interpreted, or
explained to an Italian citizen?  This problem does not arise in the
linear case when $E=1$. In this case there is a divisor-value $D =
830\,000$ with the clear interpretation that every group of $830\,000$
Union citizens accounts for one seat in Parliament (subject to
rounding).

%%%%%%%%%%%%%%%%%%%%%%%%%%%%%%%%%%%%%%%%%%%%%%%%%%%%%%%%%%%%%%%%%%%%%%%%%%%

\begin{table}[tp]
\small
\centering
\begin{tabular*}{\textwidth}{lrrrrrrrrr}
{\it Member State}&\it Population&\multicolumn{8}{c}{\it Exponent $E$ for Power-Weighted Variant}\\ 
(EU27)&(1.1.2011)&0.01&0.893&0.899&0.9&0.91&0.912&1& 27.5\\
\seprule
Germany	      & 81\,802\,257&28&95& 95&*96& 96&97&104&595\\
France	      & 64\,714\,074&28&78& 78& 78&*79&79& 83&  6\\
UK            & 62\,008\,048&28&75&*76& 76& 76&76& 80&  6\\
Italy	      & 60\,340\,328&28&74& 74& 74& 74&74& 78&  6\\
Spain	      & 45\,989\,016&28&59& 59& 59& 59&59& 61&  6\\
Poland	      & 38\,167\,329&28&51& 51& 51& 51&51& 51&  6\\
Romania	      & 21\,462\,186&28&33& 33&*32& 32&32& 31&  6\\
Netherlands   & 16\,574\,989&28&27& 27& 27& 27&27& 25&  6\\
Greece	      & 11\,305\,118&28&21& 21& 21& 21&20& 19&  6\\
Belgium	      & 10\,839\,905&28&20& 20& 20& 20&20& 19&  6\\
Portugal      & 10\,637\,713&28&20& 20& 20& 20&20& 18&  6\\
Czech Republic& 10\,506\,813&28&20& 20& 20& 20&20& 18&  6\\
Hungary	      & 10\,014\,324&28&19& 19& 19& 19&19& 18&  6\\
Sweden	      &  9\,340\,682&28&18& 18& 18& 18&18& 17&  6\\
Austria	      &  8\,375\,290&28&17& 17& 17& 17&17& 16&  6\\
Bulgaria      &  7\,563\,710&28&16& 16& 16& 16&16& 15&  6\\
Denmark	      &  5\,534\,738&28&14&*13& 13& 13&13& 12&  6\\
Slovakia      &  5\,424\,925&28&13& 13& 13& 13&13& 12&  6\\
Finland	      &  5\,351\,427&28&13& 13& 13& 13&13& 12&  6\\
Ireland	      &  4\,467\,854&28&12& 12& 12& 12&12& 11&  6\\
Lithuania     &  3\,329\,039&28&11& 11& 11&*10&10& 10&  6\\
Latvia	      &  2\,248\,374&28& 9&  9&  9&  9& 9&  8&  6\\
Slovenia      &  2\,046\,976&27& 9&  9&  9&  9& 9&  8&  6\\
Estonia	      &  1\,340\,127&27& 8&  8&  8&  8& 8&  7&  6\\
Cyprus	      &   \,803\,147&27& 7&  7&  7&  7& 7&  6&  6\\
Luxembourg    &   \,502\,066&27& 6&  6&  6&  6& 6&  6&  6\\
Malta	      &   \,412\,970&27& 6&  6&  6&  6& 6&  6&  6\\
\seprule
{\it Sum}       &501\,103\,425&751&751&751&751&751&751&751&751\\
\end{tabular*}

\bigskip
\caption{{\it Power-weighted variant, below and beyond the  $96$-seat
restriction.}  When $E = 0.01$, seats are apportioned almost equally.  When
$E=0.9$ and $E=0.91$, Germany is allocated 96 seats, stars * indicating
intervening seat transfers.  With $E=1$, population figures are unweighted.
When $E=27.5$, seats are apportioned as unequally as possible.}
\label{table1}
\end{table}
\clearpage

%%%%%%%%%%%%%%%%%%%%%%%%%%%%%%%%%%%%%%%%%%%%%%%%%%%%%%%%%%%%%%%%%%%%%%%%%%%

\begin{table}[tp]
\small
\centering
\begin{tabular*}{\textwidth}{lrrrrrrrr}
&&&&A&B& C&D& E\cr
EU27&\it Population&\it Popn$^{0.91}$&\it 
Popn$^{0.9}$&$\it CC$&\it Par.&$x(0.91)$&$x(0.9)$&\it Now\\
\seprule
DE& 81\,802\,257&15\,871\,442.9&13\,227\,834.7&96&96&96&96&99\\
FR& 64\,714\,074&12\,823\,567.3&10\,712\,698.1&85&80&79&78&74\\
UK& 62\,008\,048&12\,334\,675.7&10\,308\,684.6&81&78&76&76&73\\
IT& 60\,340\,328&12\,032\,419.8&10\,058\,816.8&79&76&74&74&73\\
ES& 45\,989\,016& 9\,397\,563.0& 7\,877\,505.3&62&61&59&59&54\\
PL& 38\,167\,329& 7\,931\,211.1& 6\,660\,741.7&52&52&51&51&51\\
RO& 21\,462\,186& 4\,697\,029.6& 3\,967\,405.2&32&33&32&32&33\\
NL& 16\,574\,989& 3\,712\,807.7& 3\,144\,183.8&26&27&27&27&26\\
EL& 11\,305\,118& 2\,621\,080.4& 2\,228\,166.2&19&20&21&21&22\\
BE& 10\,839\,905& 2\,522\,744.0& 2\,145\,472.4&19&20&20&20&22\\
PT& 10\,637\,713& 2\,479\,887.2& 2\,109\,421.9&18&19&20&20&22\\
CZ& 10\,506\,813& 2\,452\,102.5& 2\,086\,046.1&18&19&20&20&22\\
HU& 10\,014\,324& 2\,347\,284.3& 1\,997\,834.3&18&19&19&19&22\\
SE&  9\,340\,682& 2\,203\,152.3& 1\,876\,466.1&17&18&18&18&20\\
AT&  8\,375\,290& 1\,994\,940.1& 1\,700\,982.6&16&16&17&17&19\\
BG&  7\,563\,710& 1\,818\,229.6& 1\,551\,891.6&15&15&16&16&18\\
DK&  5\,534\,738& 1\,368\,416.6& 1\,171\,621.6&12&13&13&13&13\\
SK&  5\,424\,925& 1\,343\,687.6& 1\,150\,679.5&12&13&13&13&13\\
FI&  5\,351\,427& 1\,327\,111.3& 1\,136\,639.2&12&13&13&13&13\\
IE&  4\,467\,854& 1\,126\,134.0&  \,966\,249.0&11&11&12&12&12\\
LT&  3\,329\,039&  \,861\,608.9&  \,741\,458.8&10&10&10&11&12\\
LV&  2\,248\,374&  \,602\,837.7&  \,520\,812.9& 8& 8& 9& 9& 9\\
SI&  2\,046\,976&  \,553\,493.6&  \,478\,631.7& 8& 8& 9& 9& 8\\
EE&  1\,340\,127&  \,376\,446.1&  \,326\,912.4& 7& 7& 8& 8& 6\\
CY&   \,803\,147&  \,236\,245.5&  \,206\,212.8& 6& 7& 7& 7& 6\\
LU&   \,502\,066&  \,154\,060.9&  \,135\,109.1& 6& 6& 6& 6& 6\\
MT&   \,412\,970&  \,128\,969.2&  \,113\,325.2& 6& 6& 6& 6& 6\\
\seprule
\it Sum&501\,103\,425&&&751&751&751&751&754\\
\end{tabular*}

\bigskip
\caption{{\it Comparison of five EP compositions.}  
Column A is the Cambridge Compromise ``5+Upw'' of \cite{grim}.  
Column B is the parabolic allotment of \cite{ramirez}.  
Columns C and D are power-weighted variants, with exponents $0.91$ and $0.9$
respectively, and associated population indices as shown. 
Column E is the {\it status quo}.}
\label{table2}
\end{table}
\clearpage

%%%%%%%%%%%%%%%%%%%%%%%%%%%%%%%%%%%%%%%

\section{Discussion} \label{sec:4}

The legal principles for the composition of and election to the European Parliament 
(EP) may be found in Article 14 of the Treaty of Lisbon  \cite{teu}:

\begin{enumerate}

\item[(II 1)] The EP shall be composed of representatives of the Union's
citizens.  They shall not exceed seven hundred and fifty in number, plus
the President.  Representation of citizens shall be degressively
proportional, with a minimum threshold of six members per Member State.  No
Member State shall be allocated more than ninety-six seats.

\item[(III)] The members of the EP shall be elected for a term of five
years by direct universal suffrage in a free and secret ballot.
\end{enumerate}
The concept of  ``equality'' appears in several other articles  of the Treaty
(including Articles 2, 4, 9, 10), but is notably absent from Article 14.

Moreover, under the Treaty of Lisbon the task of deciding on the EP's
composition has become a matter for 
secondary (parliamentary) law, and therefore subject
to challenge in
the Court of Justice of the European Union.  

The European Commission for Democracy through Law [Venice Commission],  an institution of the Council of Europe, lists in \cite{ecd} five principles of
Europe's electoral heritage:

\begin{enumerate}
\item[(I)] The five principles underlying Europe's electoral heritage are
{\it universal, equal, free, secret and direct suffrage}.
\end{enumerate}
The 27 Member States of the Union count among the 47 members of
the Council of Europe, and as such endorse the Venice Commission's Code of
Good Practice in Electoral Matters.  Furthermore, the Union has resolved to
accede to the European Convention for the Protection of Human Rights and
Fundamental Freedom, and will then become accountable to the European Court
of Human Rights.

The concept of electoral equality  
is not self-explanatory, but requires interpretation.  
Constitutional courts take pains to distinguish between large
electorates (such as a national parliament), medium-size electorates (such as
a provincial legislature), and small electorates (as in local communities).
The European Union is on a scale in excess of these, and the notions of electoral
equality entertained within the 27 Member States cannot be extended automatically
to the entirety of the Union.  The EP has thus substantial freedom in specifying 
an interpretation of the concept of electoral equality.

The two-stage process of the Cambridge Compromise may be interpreted as
a type of ``dual electoral equality''.  This
dual concept is a merger of the {\it one state, one vote} rule of
international law, and of the {\it one person, one vote} principle of equal
representation of citizens.  There is a potential ambiguity in the term
``state'' over whether it refers to {\it government} or to {\it people}.
Our belief that many Members of the EP consider themselves 
representatives of their {\it citizenry} is supported
by the statement of Article 14 II that says just this.  We shall thus write here of
``citizenries'' rather than of ``states''.

The Cambridge Compromise  merges the two aspects of electoral
equality, and may be interpreted as follows.  
The base component is directed towards equality of citizenries,
and the proportional allocation of the remaining seats
is aimed at equality among Union citizens.  
In contrast, it is considerably harder to justify the power-weighted variant of the
Cambridge Compromise, since it violates the principle
of equal suffrage as set out by the Venice Commission:

\begin{enumerate}
\item [(2.2)] {\it Equal voting power:} seats must be evenly distributed
between the constituencies.
\end{enumerate}
There seems little doubt in the context of the EP that the word ``constituencies'' 
should be interpreted as 
``Member States".  Yet the power-weighted apportionment  
is decidedly unequal\footnote{In the terminology of the German Federal Constitutional Court, equal voting power is termed {\it Z\"ahlwertgleichheit}.}.  
For example, whereas the Italian population-index
accounts for just a sixth of its population, the Maltese index accounts for one quarter.  

A further principle of electoral affairs, the {\it continuity
  principle}, shields the legislator from abstract rules in situations
where more sensitive action is needed.  During a period of significant
institutional development, the EP may adopt the power-weighted variant
5+Pwr+Upw as a step along a continuous transition from the negotiated
{\it status quo} composition to the constitutionally principled
Cambridge Compromise.  This variant honours the equality principle for
citizenries, and converges for an enlarged Union towards the equality
principle for individual citizens.  The transitional period with $E<1$
is justifiable by an appeal to the principle of continuity.

The above discussion may be summarised as follows.
While the power-weighted variant of the Cambridge Compromise 
is technically sound and feasible, it is in
conflict with the principle of equal voting power.  Since the
power-weighting will diminish as the Union grows, this variant may
be justified by the continuity principle.

We close with a word of caution.
The lack of transparency of the power-weighted variant 
may be more harmful than helpful to public reception.
The current rapporteur Andrew Duff MEP
has proposed an intermediate step of negotiation as a bridge
between the {\it status quo}
and the Cambridge Compromise.  This
pragmatic proposal may be superior in communication and implementation.

\newpage
\bibliographystyle{amsplain}

\end{document}